\RequirePackage{fix-cm}
\documentclass[smallextended]{svjour3}       

\smartqed  

\usepackage{mathptmx}      
\usepackage{graphicx}
\usepackage{latexsym,amssymb,amsmath} 
\usepackage{mathrsfs} 

\def\abs#1{\left\vert #1 \right\vert}
\def\allpoly{\mbox{$\re\langle X \rangle$}}
\def\allpolyell{\mbox{$\re^{\ell}\langle X \rangle$}}
\def\allpolyx0degn{\mbox{$P_n$}}

\def\allseries{\mbox{$\re\langle\langle X \rangle\rangle$}}

\def\allseriesell{\mbox{$\re^{\ell} \langle\langle X \rangle\rangle$}}

\def\allseriesellLC{\mbox{$\re^{\ell}_{LC}\langle\langle X \rangle\rangle$}}
\def\allseriesellGC{\mbox{$\re^{\ell}_{GC}\langle\langle X \rangle\rangle$}}

\def\allseriesX1{\mbox{$\re [[ X_1 ]]$}}
\def\bfem#1{{\bf \em #1}} 

\def\bull{\rule{0.08in}{0.08in}} 

\newcommand{\comment}[1]{} 

\def\eqref#1{(\ref{#1})} 

\def\expup{{\rm e}}

\def\Lpm{L_{\mathfrak{p}}^m}

\def\Lpme{L^m_{\mathfrak{p},e}}

\def\norm#1{\Vert#1\Vert}

\def\openbull{\framebox[0.08in][c]{$\;$}} 

\def\re{{\mathbb R}} 

\def\sameau{\rule[0.017in]{0.2in}{0.012in}}

\def\shuffle{{\scriptscriptstyle \;\sqcup \hspace*{-0.05cm}\sqcup\;}}

\def\supp{{\rm supp}}

\def\begals{\[\begin{aligned}}
\def\endals{\end{aligned}\]}
\def\begce{\begin{center}}
\def\endce{\end{center}}
\def\begar{\begin{array}}
\def\endar{\end{array}}
\def\begeq{\begin{equation}}
\def\endeq{\end{equation}}
\def\begdi{\begin{displaymath}}
\def\enddi{\end{displaymath}}
\def\begdis{\begin{eqnarray*}}
\def\enddis{\end{eqnarray*}}
\def\begeqa{\begin{eqnarray}}
\def\endeqa{\end{eqnarray}}
\def\begdes{\begin{description}}
\def\enddes{\end{description}}
\def\begit{\begin{itemize}}
\def\endit{\end{itemize}}
\def\begen{\begin{enumerate}}
\def\enden{\end{enumerate}}
\def\beglar{\left[\begin{array}}
\def\endrar{\end{array}\right]}
\def\begle{\begin{mylemma}}
\def\endle{\end{mylemma}}
\def\begde{\begin{mydefinition}}
\def\endde{\end{mydefinition}}
\def\begth{\begin{mytheorem}}
\def\endth{\end{mytheorem}}
\def\begco{\begin{mycorollary}}
\def\endco{\end{mycorollary}}
\def\begprop{\begin{myproposition}}
\def\endprop{\end{myproposition}}
\def\begex{\begin{myexample}}
\def\endex{\hfill\openbull \end{myexample} \vspace*{0.15in}}
\def\begexer{\begin{myexercise}}
\def\endexer{\end{myexercise}}

\def\begres{\noindent{\bf Remarks}:\begin{enumerate}}
\def\endres{\end{enumerate} \par}
\def\begpr{\noindent{\em Proof:}$\;\;$}
\def\endpr{\hfill\bull \vspace*{0.15in}}
\def\begtab{\begin{tabular}}
\def\endtab{\end{tabular}}
\def\rref#1{(\ref{#1})}
\def\allseriesA{\mbox{$\re\!\ll\!A\!\gg$}}
\def\allseriesA'{\mbox{$\re\!\ll\!A'\!\gg$}}

\def\begce{\begin{center}}
\def\endce{\end{center}}
\def\begar{\begin{array}}
\def\endar{\end{array}}
\def\begeq{\begin{equation}}
\def\endeq{\end{equation}}
\def\begdi{\begin{displaymath}}
\def\enddi{\end{displaymath}}
\def\begdis{\begin{eqnarray*}}
\def\enddis{\end{eqnarray*}}
\def\begeqa{\begin{eqnarray}}
\def\endeqa{\end{eqnarray}}
\def\begdes{\begin{description}}
\def\enddes{\end{description}}
\def\begit{\begin{itemize}}
\def\endit{\end{itemize}}
\def\begen{\begin{enumerate}}
\def\enden{\end{enumerate}}
\def\beglar{\left[\begin{array}}
\def\endrar{\end{array}\right]}
\def\begle{\begin{lemma}}
\def\endle{\end{lemma}}
\def\begde{\begin{definition}}
\def\endde{\end{definition}}
\def\begth{\begin{theorem}}
\def\endth{\end{theorem}}
\def\begco{\begin{corollary}}
\def\endco{\end{corollary}}
\def\begprop{\begin{proposition}}
\def\endprop{\end{proposition}}
\def\begex{\begin{example}}
\def\endex{\hfill\openbull \end{example} \vspace*{0.1in}}
\def\begexer{\begin{exercise}}
\def\endexer{\end{exercise}}

\def\begres{\noindent{\bf Remarks}:\begin{enumerate}}
\def\endres{\end{enumerate} \par}
\def\begpr{\noindent{\em Proof:}$\;\;$}
\def\endpr{\hfill\bull \vspace*{0.1in}}
\def\begtab{\begin{tabular}}
\def\endtab{\end{tabular}}
\def\rref#1{(\ref{#1})}

\journalname{Numerische Mathematik}

\allowdisplaybreaks[4]

\begin{document}

\title{Discrete-Time Approximations of Fliess Operators\thanks{
The first author was supported by grant SEV-2011-0087 from the Severo Ochoa
Excellence Program at the Instituto de Ciencias Matem\'{a}ticas in Madrid,
Spain. The third author was supported by Ram\'on y Cajal research grant
RYC-2010-06995 from the Spanish government. This research was also supported by
a grant from the BBVA Foundation.
}
}

\author{W. Steven Gray \and
        Luis A. Duffaut Espinosa \and
        Kurusch Ebrahimi-Fard
}

\institute{W. Steven Gray \at
              Department of Electrical and Computer Engineering,
              Old Dominion University, Norfolk, Virginia 23529, USA \\
              Tel.: +1-757-683-4671\\
              Fax: +1-757-683-3220\\
              \email{sgray@odu.edu}
           \and
           Luis A. Duffaut Espinosa \at
              Department of Electrical and Computer Engineering,
              George Mason University, Fairfax, Virginia 22030, USA
           \and
           Kurusch Ebrahimi-Fard \at
              Instituto de Ciencias Matem\'{a}ticas,
              Consejo Superior de Investigaciones Cient\'{\i}ficas,
              C/ Nicol\'{a}s Cabrera, no.~13-15, 28049 Madrid, Spain
}


\maketitle

\begin{abstract}
A convenient way to represent a nonlinear input-output system in control theory is
via a Chen-Fliess functional expansion or Fliess operator.
The general goal of this paper is to describe how to approximate Fliess operators with iterated sums and to provide
accurate error estimates for two different scenarios, one where the series coefficients are growing at a local convergence rate, and the other where
they are growing at a global convergence rate. In each case, it is shown that the error estimates are achievable in the sense that worst case
inputs can be identified which hit the error bound.
The paper then focuses on the special case where the operators are
rational, i.e., they have rational generating series. It is shown in this situation that the iterated sum approximation
can be realized by a discrete-time state space model which is a rational function of the input and state affine. In addition,
this model comes from a specific discretization of the bilinear realization of the rational Fliess operator.
\end{abstract}

\section{Introduction}

A convenient way to represent a nonlinear input-output system in control theory is
via a Chen-Fliess functional expansion or Fliess operator \cite{Fliess_81,Fliess_83,Isidori_95}.
This series of weighted iterated integrals of the input functions exhibits considerable
algebraic structure that can be used, for example, to describe system
interconnections \cite{Gray-et-al_SCL14,Gray-Li_05} and to perform system inversion \cite{Gray-et-al_CDC15,Gray-et-al_Auto14}.
On the other hand, in the context
of numerical simulation and approximation, it is less clear how such a representation can be
utilized efficiently. In guidance applications, for example, piecewise constant approximations of the input have been used
in combination with a truncated version of the series to find acceptable solutions to
specific problems \cite{He-etal_2013,Yao-etal_2008}. But no a priori error estimates are provided for this approach.
Passing through a discrete-time approximation of an equivalent state space model is also an option, but not every
Fliess operator is realizable by a system of differential equations \cite{Fliess_83}.
One hint to the general problem of approximating Fliess operators was provided
by Gr\"{u}ne and Kloeden in \cite{Grune-Kloeden_01}, where it was shown that iterated integrals
can be well approximated by iterated sums. But there is a considerable jump in going from
approximating a single iterated integral to approximating an infinite sum of
such integrals. In particular, the error estimates for each iterated integral have to be
precise enough to yield an accurate error estimate for the whole operator.
Further complicating the picture is the fact that in practice only finite sums can be computed.
So an independent truncation error also has to be accounted for.

The general goal of this paper is to describe how to approximate Fliess operators with iterated sums and to provide
accurate error estimates for different scenarios. The starting point is to develop a refinement of the error estimate in
\cite[Lemma 2]{Grune-Kloeden_01} for a single iterated integral. This is done largely using Chen's Lemma \cite{Chen_52}.
After this, two specific cases are considered, one in which the series coefficients are growing at a local convergence rate, and the other where
they are growing at a global convergence rate \cite{Gray-Wang_SCL02}. Each case yields different error estimates, and several simulation examples are
given to demonstrate the results. In particular, it is shown that the error estimates are achievable in the sense that worst case
inputs can be identified which hit the error bound.
The paper then focuses on the special case where the operators are
rational, i.e., have rational generating series \cite{Berstel-Reutenauer_88}. In particular, it is shown that the iterated sum approximation
of a rational Fliess operator can be realized by a discrete-time state space model which is a rational function of the input and state affine.
This means that the approximating iterated sums do not
have to be computed explicitly but can be done implicitly via a difference equation. In which case, the truncation error can
be completely avoided. It is also shown that this difference equation approach
can be viewed in terms of a specific discretization of a continuous-time bilinear realization of the rational Fliess operator.

The paper is organized as follows. First some preliminaries on Fliess operators, Chen's Lemma, and rational series are given to set the notation and
terminology. Next the notion of a discrete-time Fliess operator is developed in Section~\ref{sec:DT-Fliess-operator}.
Then the main approximation theorems are given in Section~\ref{sec:approx-theory}. In the subsequent section, the material
concerning rational operators is presented. The conclusions of the paper are given in the final section.

\section{Preliminaries}
\label{sec:prelininaries}

A finite nonempty set of noncommuting symbols $X=\{ x_0,x_1,$ $\ldots,x_m\}$ is
called an {\em alphabet}. Each element of $X$ is called a {\em
letter}, and any finite sequence of letters from $X$,
$\eta=x_{i_1}\cdots x_{i_k}$, is called a {\em word} over $X$. The
{\em length} of $\eta$, $\abs{\eta}$, is the number of letters in
$\eta$.
The set of all words with length $k$ is denoted by
$X^k$. The set of all words including the empty word, $\emptyset$,
is designated by $X^\ast$. It forms a monoid under catenation.
The set $\eta X^\ast$ is comprised of all words with the prefix $\eta$.
Any mapping $c:X^\ast\rightarrow
\re^\ell$ is called a {\em formal power series}. The value of $c$ at
$\eta\in X^\ast$ is written as $(c,\eta)$ and called the {\em coefficient} of
$\eta$ in $c$.
Typically, $c$ is
represented as the formal sum $c=\sum_{\eta\in X^\ast}(c,\eta)\eta.$
If the {\em constant term} $(c,\emptyset)=0$ then $c$ is said to be {\em
proper}. The {\em support} of $c$, $\supp(c)$, is the set of all words having
nonzero coefficients. The collection of all formal power series over $X$ is
denoted by $\allseriesell$. The subset of polynomials is written as $\allpolyell$.
Each set forms an associative $\re$-algebra under
the catenation product and a commutative and associative $\re$-algebra under
the shuffle product, denoted here by $\shuffle$. The latter is the
$\re$-bilinear extension of the shuffle product of two words, which is defined
inductively by
\begdi
(x_i\eta)\shuffle
(x_j\xi)=x_i(\eta\shuffle(x_j\xi))+x_j((x_i\eta)\shuffle \xi)
\enddi
with $\eta\shuffle\emptyset=\emptyset\shuffle\eta=\eta$ for all $\eta,\xi\in
X^\ast$ and $x_i,x_j\in X$.

\subsection{Fliess Operators}

One can formally associate with any series $c\in\allseriesell$ a causal
$m$-input, $\ell$-output operator, $F_c$, in the following manner.
Let $\mathfrak{p}\ge 1$ and $t_0 < t_1$ be given. For a Lebesgue measurable
function $u: [t_0,t_1] \rightarrow\re^m$, define
$\norm{u}_{\mathfrak{p}}=\max\{\norm{u_i}_{\mathfrak{p}}: \ 1\le
i\le m\}$, where $\norm{u_i}_{\mathfrak{p}}$ is the usual
$L_{\mathfrak{p}}$-norm for a measurable real-valued function,
$u_i$, defined on $[t_0,t_1]$.  Let $L^m_{\mathfrak{p}}[t_0,t_1]$
denote the set of all measurable functions defined on $[t_0,t_1]$
having a finite $\norm{\cdot}_{\mathfrak{p}}$ norm and
$B_{\mathfrak{p}}^m(R)[t_0,t_1]:=\{u\in
L_{\mathfrak{p}}^m[t_0,t_1]:\norm{u}_{\mathfrak{p}}\leq R\}$. Assume
$C[t_0,t_1]$ is the subset of continuous functions in $L_{1}^m[t_0,t_1]$.
Define inductively for each $\eta\in X^{\ast}$ the map $E_\eta:
L_1^m[t_0, t_1]\rightarrow C[t_0, t_1]$ by setting $E_\emptyset[u]=1$ and
letting \[E_{x_i\bar{\eta}}[u](t,t_0) =
\int_{t_0}^tu_{i}(\tau)E_{\bar{\eta}}[u](\tau,t_0)\,d\tau, \] where
$x_i\in X$, $\bar{\eta}\in X^{\ast}$, and $u_0=1$. The
input-output operator corresponding to $c$ is the {\em Fliess operator}
\begeq
F_c[u](t) =
\sum_{\eta\in X^{\ast}} (c,\eta)\,E_\eta[u](t,t_0)
\label{eq:Fliess-operator-defined}.
\endeq
If there exist real numbers $K_c,M_c>0$ such that
\begin{equation}
\abs{(c,\eta)}\le K_c M_c^{|\eta|}|\eta|!,\;\; \forall\eta\in X^{\ast},
\label{eq:local-convergence-growth-bound}
\end{equation}
then $F_c$ constitutes a well defined mapping from
$B_{\mathfrak p}^m(R)[t_0,$ $t_0+T]$ into $B_{\mathfrak
q}^{\ell}(S)[t_0, \, t_0+T]$ provided $\bar{R}:=\max\{R,T\}<1/M_c(m+1)$,
and the numbers $\mathfrak{p},\mathfrak{q}\in[1,\infty]$ are
conjugate exponents, i.e., $1/\mathfrak{p}+1/\mathfrak{q}=1$ \cite{Gray-Wang_SCL02}.
(Here, $\abs{z}:=\max_i \abs{z_i}$ when $z\in\re^\ell$.) In this case,
the operator $F_c$ is said to be {\em locally convergent} (LC), and the set of all
series satisfying \rref{eq:local-convergence-growth-bound}
is denoted by $\allseriesellLC$.
When $c$ satisfies the more stringent growth condition
\begeq
\abs{(c,\eta)}\le K_c M_c^{|\eta|},\;\; \forall\eta\in X^{\ast}, \label{eq:global-convergence-growth-bound}
\endeq
the series \rref{eq:Fliess-operator-defined}
defines an operator from the extended space
$\Lpme (t_0)$ into $C[t_0, \infty)$,  where
\begdi
\Lpme(t_0):=
\{u:[t_0,\infty)\rightarrow \re^m:u_{[t_0,t_1]}\in \Lpm[t_0,t_1],\forall t_1 \in (t_0,\infty)\},
\enddi
and $u_{[t_0,t_1]}$ denotes the restriction of $u$ to $[t_0,t_1]$ \cite{Gray-Wang_SCL02}.
In this case, the operator is said to be {\em globally convergent} (GC), and the set of all
series satisfying \rref{eq:global-convergence-growth-bound}
is designated by
$\allseriesellGC$.

\subsection{Chen's Lemma}

\begin{figure}[t]
\includegraphics[scale=0.32,angle=-90]{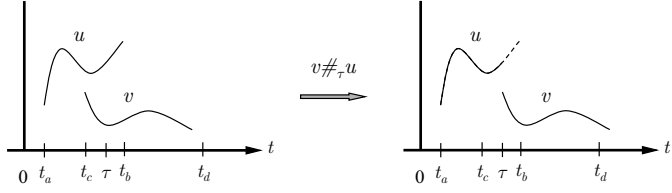}
\caption{The catenation of two inputs $u$ and $v$ at $t=\tau$.}
\label{fig:function-catenation}
\end{figure}

For a fixed $u$
consider a series in $\allseries$ of the form $P[u]=\sum_{\eta\in X^\ast} \eta E_{\eta}[u]$, which
is often referred to as a {\em Chen series}.
Given two functions $(u,v)\in L_1^m[t_a,t_b]\times L_1^m[t_c,t_d]$, their
{\em durations} are taken to be $t_b-t_a\geq 0$ and
$t_d-t_c\geq 0$, respectively, and the functions are
not defined outside their corresponding intervals. The
{\em catenation} of $u$ and $v$
at $\tau\in[t_a,t_b]$ is understood to be
\begdi
(v\#_{\tau}u)(t)=\left\{\begar{ccl}
           u(t) &:& t_a \leq t \leq \tau \\
           v((t-\tau)+t_c) &:& \tau < t \leq \tau+(t_d-t_c)
         \endar\right.
\enddi
(see Figure~\ref{fig:function-catenation}).
It is easily verified that
$
L_{1,e}^m(0)
$
is a monoid under the catenation operator.
The identity element in this case is denoted by $\mathbf{0}$ and is
equivalent to the set of functions having exactly zero duration.
The following lemma is due to Chen \cite{Chen_52}.

\begle (Chen's Lemma)
If $(u,v)\in L_1^m[0,T_1]\times L_1^m[0,T_2]$ and
$(t_1,t_2)\in[0,T_1]\times[0,T_2]$ then
\begdi
P[v](t_2)P[u](t_1)=P[v\#_{t_1} u](t_2+t_1). \\[0.1in]
\enddi
\endle

So in essence $P:L_{1,e}^m(0)\rightarrow \allseries$ acts as a monoid morphism,
where $\allseries$ is viewed as a monoid under the catenation product.

\subsection{Rational Formal Power Series}

A brief summary of rational and recognizable formal power series is useful.
The treatment here is based largely on \cite{Berstel-Reutenauer_88}.

A series $c\in\allseries$ is called {\em invertible} if there exists
a series $c^{-1}\in\allseries$ such that $cc^{-1}=c^{-1}c=1$.\footnote{The polynomial $1\emptyset$ is abbreviated throughout as $1$.} In the
event that $c$ is not proper, it is always possible to write
\begdi
c=(c,\emptyset)(1-c^{\prime}),
\enddi
where $(c,\emptyset)$ is nonzero, and $c^{\prime}\in\allseries$ is proper.
It then follows that
\begdi
c^{-1}=\frac{1}{(c,\emptyset)}(1-c^{\prime})^{-1}
=\frac{1}{(c,\emptyset)} (c^\prime)^{\ast},
\enddi
where
\begdi
(c^{\prime})^{\ast}:=\sum_{i=0}^{\infty} (c^{\prime})^i.
\enddi
In fact, $c$ is invertible if and {\em only if} $c$ is not proper.
Now let $S$ be a subalgebra of the $\re$-algebra $\allseries$
with the catenation product.
$S$ is said to be {\em rationally closed} when every invertible
$c\in S$ has $c^{-1}\in S$
(or equivalently, every proper $c^{\prime}\in S$ has $(c^{\prime})^{\ast}\in S$).
The {\em rational closure} of any
subset $E\subset\allseries$ is the smallest rationally closed
subalgebra of $\allseries$ containing $E$.

\begde
A series $c\in\allseries$ is \bfem{rational} if it belongs to the
rational closure of $\allpoly$.
\endde

It turns out that an entirely different characterization of a rational series is possible
using the following concept.

\begde \label{def:linear-representation}
A \bfem {linear} \bfem{representation}
of a series $c\in\allseries$ is any triple
$(\mu,\gamma,\lambda)$, where
\begdi
\mu:X^{\ast}\rightarrow \re^{n\times n}
\enddi
is a monoid morphism, and $\gamma,\lambda^T\in\re^{n\times 1}$ are such that
\begdi
(c,\eta)=\lambda\mu(\eta)\gamma, \;\forall\eta\in X^{\ast}.
\enddi
The integer $n$ is the dimension of the representation.
\endde

\begde
A series $c\in\allseries$ is called \bfem{recognizable} if
it has a linear representation.
\endde

\begth {\rm (Sch\"{u}tzenberger)} \label{th:Schutzenberger}
A formal power series is rational if and only if it is recognizable.
\endth

The next concept
provides an explicit way of constructing a linear representation of a rational series.
Define for any $x_i\in X$, the left-shift operator, $x_i^{-1}(\cdot)$, on $X^\ast$ by
$x_i^{-1}(x_i\eta)=\eta$ with $\eta\in X^\ast$ and zero otherwise.
Higher order shifts are defined inductively via $(x_i\xi)^{-1}(\cdot)=\xi^{-1}x_i^{-1}(\cdot)$,
where $\xi\in X^\ast$. The left-shift operator is assumed to act linearly on $\allseries$.

\begde
A subset $V\subset\allseries$ is called \bfem{stable} when $\xi^{-1}(c)\in V$
for all $c\in V$ and $\xi\in X^{\ast}$.
\endde

\begth \label{th:stable-rational-subspace}
A series $c\in\allseries$ is rational/recognizable if and only if
there exists a stable finite dimensional $\re$-vector subspace of
$\allseries$ containing $c$.
\endth

\section{Discrete-Time Fliess Operators}
\label{sec:DT-Fliess-operator}

Let $u\in L_1^m[0,T]$ for some finite $T>0$.
Following \cite{Grune-Kloeden_01},
select some integer $L\geq 1$ and with $\Delta:=T/L$ define the sequence
\begin{equation} \label{eq:uhat-def}
\hat{u}_i(N)=\int_{(N-1)\Delta}^{N\Delta} u_i(t)\,dt,\;\;i=0,1,\ldots,m
\end{equation}
where $N\in[1,L]$.
Observe in particular that $\hat{u}_0(N)=\Delta$ since $u_0=1$.
The corresponding iterated sum  for any $x_i\in X$ and $\eta\in X^\ast$ is defined inductively by
\begdi \label{eq:interated-sum-approx-for-iterated-integral}
S_{x_i\eta}[\hat{u}](N)=\sum_{k=1}^N \hat{u}_{i}(k)S_{\eta}[\hat{u}](k)
\enddi
with $S_\emptyset[\hat{u}](N):=1$.
The following lemma gives an alternative description of $S_\eta$ which will be useful later.

\begle \label{le:S-as-Cauchy-type-sum} For any $N\in[1,L]$ and $\eta\in X^\ast$
\begdi
S_\eta[\hat{u}](N)=\Delta^{\abs{\eta}}\sum_{\xi_N\cdots\xi_1=\eta} u_{\xi_N}(N)\cdots u_{\xi_1}(1),
\enddi
where $u_i(k):=\hat{u}_i(k)/\Delta$,
$u_{x_{i_1}\cdots x_{i_r}}(k):=u_{i_1}(k)\cdots u_{i_r}(k)$, $u_\emptyset(k):=1$,
and the summation is over all partitions of $\eta$ having N subwords $\xi_{k}\in X^\ast$
(so some subwords can be empty).
\endle

\begpr
The proof is by induction on the length of $\eta$. For the empty word the equality holds trivially.
When $\eta=x_i$ observe that
\begdi
S_{x_i}[\hat{u}](N)=\sum_{k=1}^N \hat{u}_{i}(k)=\Delta \sum_{k=1}^N u_{i}(k)
= \Delta^{} \sum_{\xi_N\cdots\xi_1=x_i} u_{\xi_N}(N)\cdots u_{\xi_1}(1).
\enddi
Now assume the claim holds for all words up to length $j\geq 0$.
If $\eta\in X^j$ then
\begin{align*}
S_{x_i\eta}[\hat{u}](N)&=\sum_{k=1}^N \hat{u}_{i}(k)S_{\eta}[\hat{u}](k)
=\sum_{k=1}^N \Delta u_{i}(k) \Delta^j
\sum_{\xi_k\cdots\xi_1=\eta} u_{\xi_k}(k)\cdots u_{\xi_1}(1) \\
&=\Delta^{j+1}
\sum_{\xi_N\cdots\xi_1=x_i\eta} u_{\xi_N}(N)\cdots u_{\xi_1}(1),
\end{align*}
which proves the lemma.
\endpr

The next definition provides the main class of discrete-time approximators used throughout the paper.
In the most general context,
the set of admissible inputs will be drawn from the real sequence space
$l_\infty^{m+1}[N_0]:=\{\hat{u}=(\hat{u}(N_0),\hat{u}(N_0+1),\ldots): \abs{\hat{u}(N)}<\hat{R}_{\hat{u}}<\infty,\,\forall N\geq N_0\}$,
where $\abs{\hat{u}(N)}:=\max_{i=0,1,\ldots,m}\abs{\hat{u}_i(N)}$.
In which case, $\norm{\hat{u}}_\infty:=\sup_{N\geq N_0}\abs{\hat{u}(N)}$ is always finite.
To be consistent with \rref{eq:uhat-def}, it is assumed throughout that $\hat{u}_0$ is a constant input.
Define a ball of radius $\hat{R}$ in $l_\infty^{m+1}[N_0]$ as $B_\infty^{m+1}[N_0]$ $(\hat{R})=\{\hat{u}\in l_\infty^{m+1}[N_0]: \norm{\hat{u}}_\infty\leq \hat{R} \}$.
The subset of finite sequences over $[N_0,N_f]$ is denoted by $B_\infty^{m+1}[N_0,N_f](\hat{R})$.

\begde
For any $c\in\allseriesell$, the corresponding {\bf discrete-time Fliess operator} defined on $l_\infty^{m+1}[1]$ is
\begeq \label{eq:DT-Fliess-operator-defined}
\hat{y}(N)=\hat{F}_c[\hat{u}](N)=\sum_{\eta\in X^\ast} (c,\eta)S_{\eta}[\hat{u}](N).
\endeq
\endde

Before considering the approximation problem, it is necessary to
introduce
various sufficient conditions for convergence of such operators.
The following lemma is essential.

\begle \label{le:Seta-bound}
If $\hat{u}\in B_\infty^{m+1}[1](\hat{R})$ then for any $\eta\in X^\ast$
$$
\abs{S_\eta[\hat{u}](N)}\leq \hat{R}^{\abs{\eta}} {N-1+\abs{\eta} \choose \abs{\eta}}. \\[0.1in]
$$
\endle

\begpr
If $\eta=x_{i_j}\cdots x_{i_1}$ then observe for any $N\geq 1$
\begin{align*}
\abs{S_\eta[\hat{u}](N)}
&= \abs{\sum_{k_j=1}^N \hat{u}_{i_j}(k_j)\sum_{k_{j-1}=1}^{k_j} \hat{u}_{i_{j-1}}(k_{j-1})\cdots \sum_{k_1=1}^{k_2} \hat{u}_{i_1}(k_1)} \\
&\leq \sum_{k_j=1}^N \abs{\hat{u}_{i_j}(k_j)}\sum_{k_{j-1}=1}^{k_j} \abs{\hat{u}_{i_{j-1}}(k_{j-1})}\cdots \sum_{k_1=1}^{k_2} \abs{\hat{u}_{i_1}(k_1)} \\
&\leq \hat{R}^{\abs{\eta}} \sum_{k_j=1}^N \sum_{k_{j-1}=1}^{k_j} \cdots \sum_{k_1=1}^{k_2} 1 \\
&= \hat{R}^{\abs{\eta}} {N-1+\abs{\eta} \choose \abs{\eta}},
\end{align*}
using the fact that the final nested sum above has ${N-1+\abs{\eta} \choose \abs{\eta}}$ terms \cite{Butler-Karasik_10}.
\endpr

Since the upper bound on $\abs{S_\eta[\hat{u}](N)}$ in this lemma is achievable on $B_\infty^{m+1}[1](\hat{R})$,
it is not difficult to see that when the generating
series $c$ satisfies the growth bound \rref{eq:local-convergence-growth-bound}, the series \rref{eq:DT-Fliess-operator-defined} defining $\hat{F}_c$ can
diverge. For example, if $(c,\eta)=K_cM_c^{\abs{\eta}}\abs{\eta}!$ for all $\eta\in X^\ast$, and $\hat{u}$ is such a maximizing input then
\begin{align*}
\hat{F}[\hat{u}](N)
&= K_c\sum_{\eta\in X^\ast}M_c^{\abs{\eta}}\abs{\eta}! \hat{R}^{\abs{\eta}}{N-1+\abs{\eta} \choose \abs{\eta}} \\
&=K_c\sum_{j=0}^\infty (M_c(m+1)\hat{R})^j ((N-1+j)\cdots (N+1)N),
\end{align*}
which will diverge even if $M_c(m+1)\hat{R} <1$. The next theorem shows that this problem is averted when $c$ satisfies the stronger
growth condition \rref{eq:global-convergence-growth-bound}.

\begth
Suppose $c\in\allseriesell$ has coefficients which satisfy \rref{eq:global-convergence-growth-bound}.
Then there exists a real number $\hat{R}>0$ and an integer $L\geq 1$ such that for
each $\hat{u}\in B_\infty^{m+1}[1,L](\hat{R})$, the series \rref{eq:DT-Fliess-operator-defined}
converges absolutely and uniformly on $[1,L]$.
\endth

\begpr
Fix $L\geq 1$ and select any $N\in[1,L]$.
In light of Lemma~\ref{le:Seta-bound}, if $\abs{\eta}\gg N$ then
$$
\abs{S_\eta[\hat{u}](N)}
\lesssim \hat{R}^{\abs{\eta}}\frac{1}{(N-1)!}.
$$
From the assumed coefficient bound it follows that
\begin{align*}
\abs{\hat{F}_c(\hat{u})(N)}&\leq \sum_{j=0}^\infty \sum_{\eta\in X^j} \abs{(c,\eta)} \abs{S_\eta[\hat{u}](N)}
\lesssim\sum_{j=0}^\infty K_c(M_c(m+1))^j\, \hat{R}^j \frac{1}{(N-1)!} \\
&=\frac{1}{(N-1)!}\frac{K_c}{1-M_c(m+1)\hat{R}},
\end{align*}
provided $\hat{R}<1/M_c(m+1)$.
Since $\hat{u}_0$ is constant on $[1,L]$, an upper bound on $L$ is also implied.
\endpr

The final convergence theorem shows that the restriction on the norm of $\hat{u}$ can be removed if an even more stringent
growth condition is imposed on $c$.

\begth
Suppose $c\in\allseriesell$ has coefficients which satisfy
\begdi \label{eq:exp-growth-bound}
\abs{(c,\eta)}\leq K_cM_c^{\abs{\eta}}\frac{1}{\abs{\eta}!},\;\;\eta\in X^\ast
\enddi
for some real numbers $K_c,M_c>0$.
Then for every $\hat{u}\in l_\infty^{m+1}[1]$, the series \rref{eq:DT-Fliess-operator-defined}
converges absolutely and uniformly on $[1,\infty)$.
\endth

\begpr
Following the same argument as in the proof of the previous theorem, it is clear
for any $\hat{u}\in l_\infty^{m+1}[1]$ and $N\geq 1$ that
\begdi
\abs{\hat{F}_c(\hat{u})(N)}
\lesssim \sum_{j=0}^\infty K_c(M_c(m+1))^j \frac{1}{j!}\, \|\hat{u}\|_\infty^j \frac{1}{(N-1)!}
=\frac{K_c}{(N-1)!}\expup^{M_c(m+1)\|\hat{u}\|_\infty}.
\enddi
\endpr

Assuming the
analogous definitions for local and global convergence of the operator $\hat{F}_c$,
note the incongruence between the convergence conditions for continuous-time and discrete-time Fliess operators as
summarized in Table~\ref{tbl:convergence-cond-summary}. In each case, for a fixed $c$, the sense in which $\hat{F}_c$ converges is {\em weaker} than that
for $F_c$.
This is not entirely surprising given that the input $\hat{u}$ in the approximation setting is viewed as the increments of the integral of
$u$ rather than $u$ itself. But the real source of this dichotomy is the observation in Lemma~\ref{le:Seta-bound} that iterated sums of $\hat{u}$ do not grow as a function of word length
like $1/\abs{\eta}!$, which is the case for iterated integrals. As shown in the next section, however, this difference
in convergence behavior does not provide any serious
impediment to using discrete-time Fliess operators
as approximators for their continuous-time counterparts.

\begin{table}[h]
\caption{Summary of convergence conditions for $F_c$ and $\hat{F}_c$.}
\label{tbl:convergence-cond-summary}
\renewcommand{\arraystretch}{2}
\begin{center}
\begin{tabular}{|c|c|c|} \hline
 growth rate & $F_c$ & $\hat{F}_c$ \\ \hline
$\abs{(c,\eta)}\le K_c M_c^{|\eta|}|\eta|!$ & LC & divergent \\ \hline
$\abs{(c,\eta)}\le K_c M_c^{|\eta|}$ & GC & LC \\ \hline
$\abs{(c,\eta)}\le K_c M_c^{|\eta|}\frac{1}{|\eta|!}$ & at least GC & GC \\[0.03in] \hline
\end{tabular}
\end{center}
\end{table}

\section{Approximating Fliess Operators}
\label{sec:approx-theory}

\subsection{Iterated Integrals}

The starting point for the approximation theory is the observation that
$E_{x_i}[u](T,0)=S_{x_i}[\hat{u}](L)$
for all $x_i\in X$ and the assertion of Gr\"{u}ne and Kloeden that
for any $\eta\in X^\ast$ with $|\eta|\geq 2$
\begdi \label{eq:GK-iterated-integral-approx}
S_\eta[\hat{u}](L)=E_{\eta}[u](T,0)+O\left(\frac{T^{|\eta|}}{L}\right)
\enddi
\cite[Lemma 2]{Grune-Kloeden_01}.
The following theorem gives an explicit error bound along these lines.

\begth \label{th:GK-iterated-integral-approx-improved}
Let $u\in L_1^m[0,T]$ for some finite $T>0$.
Select integer $L\geq 1$, set $\Delta:=T/L$, and define the sequence
$\hat{u}$ as in \rref{eq:uhat-def}.
For any $\eta\in X^\ast$ it follows that if $L\gg \abs{\eta}\geq 2$ then
\begdi \label{eq:GK-iterated-integral-approx-fixed}
\abs{S_\eta[\hat{u}](L)-E_\eta[u](T,0)}\lesssim
\frac{T^{\abs{\eta}}}{L}\frac{\norm{\hat{u}/\Delta}_\infty^{\abs{\eta}}}{2(\abs{\eta}-2)!}. \\[0.1in]
\enddi
\endth

\begpr
Since the input sequence $\hat{u}$ is
computed exactly from the integration of $u$,
there is no loss of generality in the computation of $S_\eta[\hat{u}](L)$ if one assumes a priori that $u$ is a piecewise constant input taking values
$u_i(t):=\hat{u}_i(N)/\Delta$ when $t\in[(N-1)\Delta,N\Delta)$ for $i=0,1,\ldots,m$.
In addition, it was shown in \cite[Lemma 2.1]{Gray-Wang_SCL02} for any $u\in L_1[0,T]$ that
\begeq \label{eq:Eeta-upper-bound-Gray-Wang}
\abs{E_\eta[u](N\Delta,(N-1)\Delta)}\leq \frac{U_0^{\abs{\eta}_{x_0}}\cdots {U_m^{\abs{\eta}_{x_m}}}}{\abs{\eta}_{x_0}!\cdots \abs{\eta}_{x_m}!},
\endeq
where $U_i:=\int_{(N-1)\Delta}^{N\Delta} \abs{u_i(\tau)}\,d\tau$, and $\abs{\eta}_{x_k}$ denotes the number of times the letter $x_k$ appears in $\eta$.
This upper bound is achieved when each $u_i$ is
constant over $[(N-1)\Delta,N\Delta)$.
Thus, the worst case error between $E_{\eta}[u](T)$ and $S_\eta[\hat{u}](L)$ occurs for piecewise constant inputs.
Applying Chen's Lemma specifically to the piecewise constant input $u=u(L)\#_{(L-1)\Delta}u(L-1)\#_{(L-2)\Delta}\cdots \#_{\Delta}u(1)$ with
$u(N):=\hat{u}(N)/\Delta$, $N=1,2,\ldots,L$, gives directly
\begin{align*}
E_\eta[u](T,0)
&=(P[u](L\Delta),\eta)
=(P[u(L)](\Delta)\cdots P[u(1)](\Delta),\eta) \\
&=\sum_{\xi_L\cdots\xi_1=\eta} E_{\xi_L}[u(L)](L\Delta,(L-1)\Delta)\cdots
E_{\xi_1}[u(1)](\Delta,0).
\end{align*}
But for any $\xi=x_{i_1}\cdots x_{i_r}$
\begdi
E_{\xi}[u(N)](N\Delta,(N-1)\Delta)=u_{i_1}(N)\cdots u_{i_r}(N)\frac{\Delta^r}{r!},
\enddi
and therefore,
\begdi
E_\eta[u](T,0)=\Delta^{\abs{\eta}}\sum_{\xi_L\cdots\xi_1=\eta} \frac{1}{\abs{\xi_L}!\cdots \abs{\xi_1}!} u_{\xi_L}(L)\cdots u_{\xi_1}(1).
\enddi
Put another way, each $P[u(N)](\Delta)$ is an exponential Lie series, so from the Baker-Campbell-Hausdorff formula the
same is true of $P[u](L\Delta)$, and $E_{\eta}[u]=(P[u],\eta)$ is a truncated version of this series.
Comparing the expression above to that for $S_\eta[\hat{u}](L)$ from Lemma~\ref{le:S-as-Cauchy-type-sum}, it follows
that if $L\gg j:=\abs{\eta}$ then
\begin{align*}
\abs{S_\eta[\hat{u}](L)-E_\eta[u](T,0)}
&\leq \Delta^j\sum_{\xi_L\cdots\xi_1=\eta} \left[1-\frac{1}{\abs{\xi_L}!\cdots \abs{\xi_1}!}\right] \abs{u_{\xi_L}(L)\cdots u_{\xi_1}(1)}\\
&\leq \norm{\hat{u}}_\infty^{j}\left(\left[\sum_{\xi_L\cdots\xi_1=\eta} 1\right] -
\left[\sum_{\xi_L\cdots\xi_1=\eta} \frac{1}{\abs{\xi_L}!\cdots \abs{\xi_1}!}\right] \right) \\
&= \norm{\hat{u}}_\infty^{j}\left({L+j-1 \choose j}-\frac{L^j}{j!}\right) \\
&=\frac{\norm{\hat{u}}_\infty^{j}}{j!}\left((L+j-1)\cdots (L+1)L-L^j\right) \\
&=\frac{\norm{\hat{u}}_\infty^{j}}{j!}\left(\frac{j(j-1)}{2}L^{j-1}+\cdots +(j-1)!L\right) \\
&\approx \frac{T^j}{L}\frac{\norm{\hat{u}/\Delta}_\infty^{j}}{2(j-2)!},
\end{align*}
which proves the lemma.
\endpr

\subsection{Locally Convergent $F_c$ }

When $c$ is locally convergent, it was shown in the previous section that $\hat{F}_c$ can diverge.
Therefore, a truncated version
of $\hat{F}_c$,
$$
\hat{F}^J_c[\hat{u}](N):=\sum_{j=0}^J \sum_{\eta\in X^j} (c,\eta) S_\eta[\hat{u}](N),
$$
is considered.
The following theorem states that the error in approximating $F_c[u](T)$ with $\hat{F}^J_c[\hat{u}](L)$
can be bounded by the sum of two errors, namely, $\hat{e}(J)$, which bounds the approximation error
between iterated integrals and iterated sums, and $e(J)$, which bounds the tail of the series defining $F_c[u](T)$,
i.e., the truncation error.

\begth \label{th:DT-FO-appproximation-LC}
Let $c\in\allseriesellLC$ with growth constants $K_c,M_c>0$.
If $u\in B_1^m(R)[0,T]$ with $\bar{R}:=\max\{R,T\}<1/M_c(m+1)$ and $L\gg J$ then
\begdi
\abs{F_c[u](T)-\hat{F}_c^J[\hat{u}](L)}\lesssim\hat{e}(J)+e(J),\\[0.1in]
\enddi
where
\begin{align*}
\hat{e}(J)&=\frac{K_c}{L} \left[\frac{\hat{s}^2}{(1-\hat{s})^3} - \frac{2J(J+1)\hat{s}^{(J+1)}}{1-\hat{s}}-\frac{J\hat{s}^{(J+2)}}{(1-\hat{s})^2}-
\frac{\hat{s}^{J+2}}{(1-\hat{s})^3}\right] \\
e(J)&=K_c\frac{s^{J+1}}{1-s}
\end{align*}
with $\hat{s}:=M_c(m+1)L\norm{\hat{u}}_\infty$ and $s:=M_c(m+1)\bar{R}$.
\endth

\begpr
Applying Theorem~\ref{th:GK-iterated-integral-approx-improved} and the assumption that $s<1$ give the following:
\begin{align*}
\abs{F_c[u](T)-\hat{F}_c^J[\hat{u}](L)}
&=\abs{\sum_{j=0}^\infty\sum_{\eta\in X^j}(c,\eta)E_\eta[u](T,0)-\sum_{j=0}^J \sum_{\eta\in X^j} (c,\eta) S_{\eta}[\hat{u}](L)} \\
&\leq\sum_{j=0}^J\sum_{\eta\in X^j}\abs{(c,\eta)}\abs{E_\eta[u](T,0)-S_{\eta}[\hat{u}](L)}+ \\
&\hspace*{0.28in}\sum_{j=J+1}^\infty \sum_{\eta\in X^j} \abs{(c,\eta)}\abs{E_\eta[u](T,0)} \\
&\lesssim\sum_{j=2}^J K_cM_c^j(m+1)^j j!\, \frac{T^j}{L}\frac{\norm{\hat{u}/\Delta}_\infty^j}{2(j-2)!}+ \\
&\hspace*{0.28in}\sum_{j=J+1}^\infty K_cM_c^j(m+1)^j j!\,\frac{\bar{R}^j}{j!} \\
&=\frac{K_c}{2L}\sum_{j=0}^J (M_c(m+1)L\norm{\hat{u}}_\infty)^j j(j-1)+ \\
&\hspace*{0.28in}K_c\sum_{j=J+1}^\infty (M_c(m+1)\bar{R})^j \\
&=\frac{K_c}{2L} \left[\frac{2\hat{s}^2}{(1-\hat{s})^3} - \frac{J(J+1)\hat{s}^{(J+1)}}{1-\hat{s}}-\frac{2J\hat{s}^{(J+2)}}{(1-\hat{s})^2}-\right. \\
&\hspace*{0.28in}\left.\frac{2\hat{s}^{J+2}}{(1-\hat{s})^3}\right]+K_c\frac{s^{J+1}}{1-s} \\
&=\hat{e}(J)+e(J),
\end{align*}
where standard formulas have been used to give closed-forms for the final two series.

\endpr

Simple examples show that it is possible to have $\hat{s}\leq s$ and $\hat{s}\geq s$, so the assumed bound $s<1$ in Theorem~\ref{th:DT-FO-appproximation-LC}
does not imply that the same holds for $\hat{s}$. But
in the event that $\hat{s}<1$  and $L\gg J\gg 1$, the following corollary gives a simplified upper bound on the approximation error.

\begco \label{co:DT-FO-appproximation-LC}
Let $c\in\allseriesellLC$ with growth constants $K_c,M_c>0$.
If $u\in B_1^m(R)[0,T]$ with $\bar{R}=\max\{R,T\}<1/M_c(m+1)$, $L\norm{\hat{u}}_\infty<1/M_c(m+1)$, and
$L\gg J\gg 1$ then
\begdi
\abs{F_c[u](T)-\hat{F}_c^J[\hat{u}](L)}\lesssim \frac{K_c (M_c(m+1)L\norm{\hat{u}}_\infty)^2}{L(1-M_c(m+1)L\norm{\hat{u}}_\infty)^3}.\\[0.1in]
\enddi
\endco

\begpr
Since $\hat{s}<1$ and $J\gg 1$ then $\hat{e}(J)\approx K_c\hat{s}^2/L(1-\hat{s})^3$. In addition, since
$s<1$ and $J\gg 1$, $e(J)\approx 0$.
\endpr

\begin{figure}[t]
\begin{center}
\includegraphics[scale=0.5]{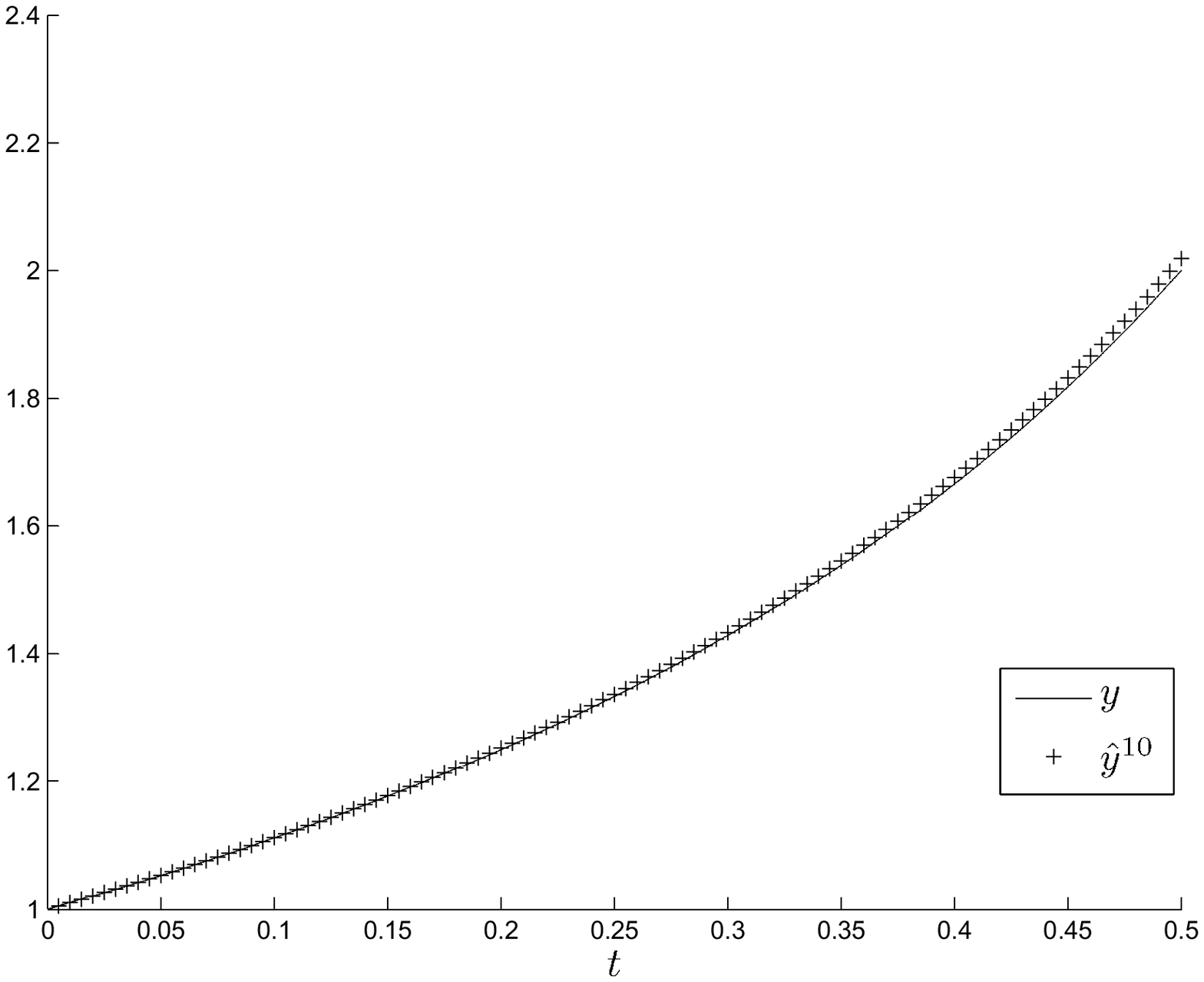}
\caption{Simulation comparing $y=F_c[1]$ to its approximation $\hat{y}^{10}=\hat{F}_c^{10}[\Delta]$ in
Example~\ref{ex:non-rational-Ferfera-DT-FO-approximation}, case 3.}
\label{fig:CTvsDT-non-rational-Ferfera-systems}
\end{center}
\end{figure}
\begin{figure}[t]
\begin{center}
\includegraphics[scale=0.5]{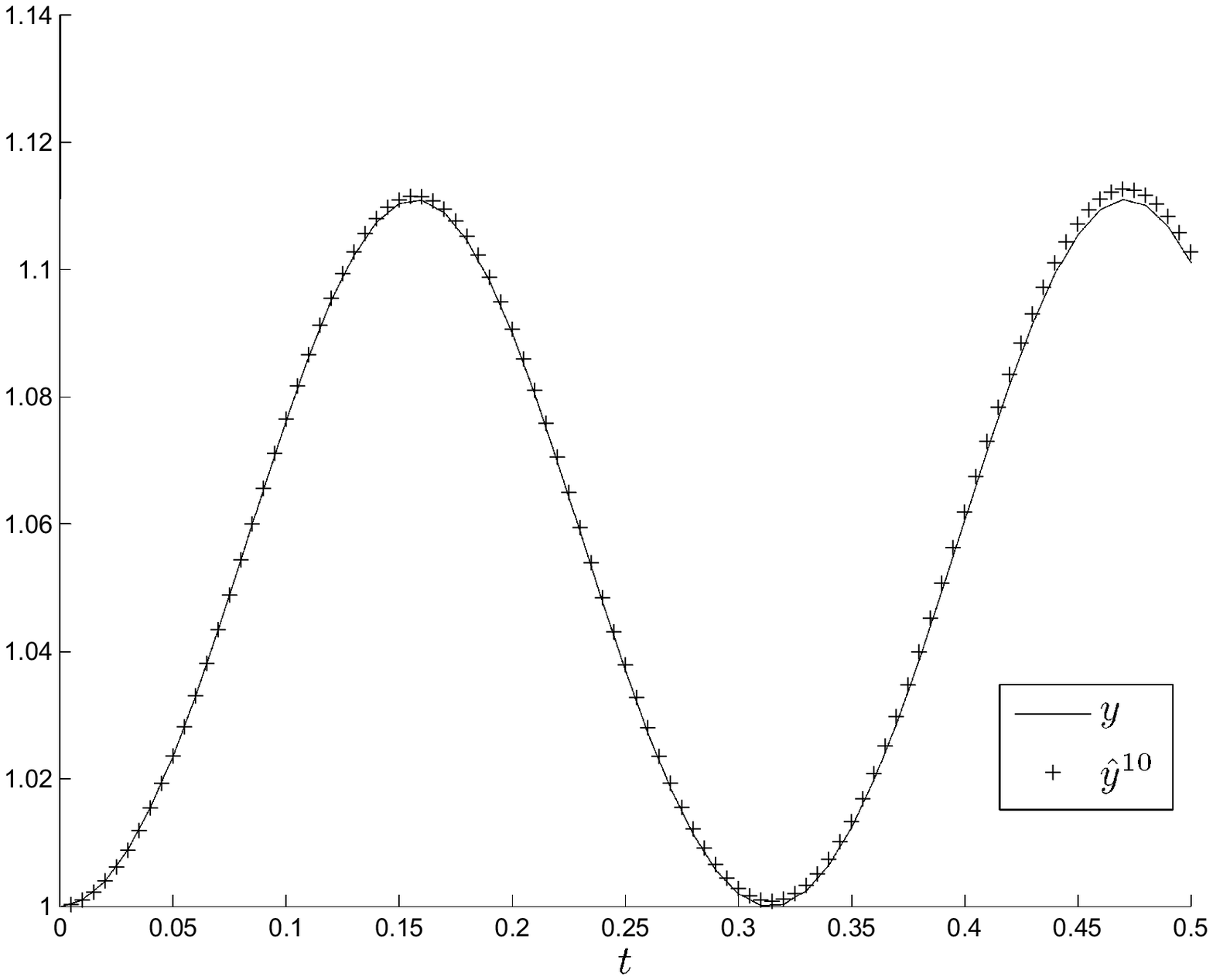}
\caption{Simulation comparing $y(t)=F_c[\sin(20t)]$ to its approximation $\hat{y}^{10}=\hat{F}_c^{10}[\hat{u}]$ in
Example~\ref{ex:non-rational-Ferfera-DT-FO-approximation}, case 6.}
\label{fig:CTvsDT-non-rational-Ferfera-systems2}
\end{center}
\end{figure}

\begex \label{ex:non-rational-Ferfera-DT-FO-approximation}
Consider the locally convergent series $c=\sum_{k\geq 0} k!\,x_1^k$ so that $K_c=M_c=1$.
Effectively, $m=0$ since $c$ only involves one letter.
It is easy to verify that $y=F_c[u]$ has the state space realization
\begdi \label{eq:CT-LC-nonrational-Ferfera-system-realization}
\dot{z}=u,\;\; z(0)=0,\;\; y=1/(1-z)
\enddi
when $\bar{R}=\max\{\norm{u}_1,T\}<1$. For example,
direct substitution for $z$ into the output equation gives
\begdi
y(t)=\sum_{j=0}^{\infty} E^j_{x_1}[u](t)= \sum_{j=0}^{\infty} E_{x_1^{\shuffle j}}[u](t)
=\sum_{j=0}^{\infty}j!\,E_{x_1^j}[u](t) = F_c[u](t).
\enddi
If the constant input $u=1$ is applied over the interval $[0,T]$ with $T<1$
then $y(T)=1/(1-T)$. On the other hand,
the discrete-time approximation $\hat{y}^J(N):=\hat{F}_c^J[\hat{u}](N)$ with $\hat{u}=\Delta$ and $N=L$ is
\begin{align*}
\hat{F}^J_c[\Delta](L)&=\sum_{j=0}^J j!\,S_{x_1^j}[\Delta](L)
= \sum_{j=0}^J j!\Delta^j\sum_{k_1+k_2+\cdots +k_L=j} \!\!\!1 \\
&=\sum_{j=0}^J j!\Delta^j  {L+j-1 \choose j} \\
&=\sum_{j=0}^J \Delta^j \left(L^j+\frac{j(j-1)}{2}L^{j-1}+\cdots + (j-1)!L\right) \\
&\approx\sum_{j=0}^J T^j + \frac{1}{2L}\sum_{j=0}^J j(j-1)T^{j} \\
&=[F_c[1](T)-e(J)]+\hat{e}(J),
\end{align*}
which is consistent with Theorem~\ref{th:DT-FO-appproximation-LC} and represents the worst case
in the sense that the upper bound \rref{eq:Eeta-upper-bound-Gray-Wang} on each iterated integral is attained.
The outputs $y$ and $\hat{y}^J$ were computed numerically over the interval $[0,0.5]$ for various choices of
$u$, $L$, and $J$. This data is summarized in Table~\ref{tbl:simulation-summary-nonrational-case} (see the last page), and the corresponding
plots for cases 3 and 6 are shown in
Figures~\ref{fig:CTvsDT-non-rational-Ferfera-systems} and
\ref{fig:CTvsDT-non-rational-Ferfera-systems2}, respectively. For this example, most of the error in
the approximation is due to the term $\hat{e}(J)$. As expected, the constant input case yields an error that is
approximately upper bounded by $\hat{e}(J)+e(J)$, while for the sinusoidal input this bound is conservative.
\endex

\subsection{Globally Convergent $F_c$}

When $c$ is globally convergent, the divergence problem for $\hat{F}_c$ is avoided provided
$\hat{u}$ is sufficiently small. But in most cases it is usually not possible to compute the infinite sum
defining $\hat{F}_c$, so once again the truncated approximator $\hat{F}_c^J$ will be utilized.
The main theorem of this section is given below. It provides an upper bound on the approximation error
in terms of the (upper) incomplete gamma function, $\Gamma(a,b):=\int_b^\infty t^{a-1}\expup^{-t}\,dt/\Gamma(a)$.


\begth \label{th:DT-FO-appproximation-GC}
Let $c\in\allseriesellGC$ with growth constants $K_c,M_c>0$.
If $u\in B_1^m(R)[0,T]$ and $L\gg J$ then
\begdi
\abs{F_c[u](T)-\hat{F}_c^J[\hat{u}](L)}\lesssim\hat{e}(J)+e(J),\\[0.1in]
\enddi
where
\begdi
\hat{e}(J)=\frac{K_c}{2L}\expup^{\hat{s}}\hat{s}^2 \Gamma(J+1,\hat{s}),\;\;\;
e(J)=K_c\expup^s(1-\Gamma(J+1,s))
\enddi
with $\hat{s}:=M_c(m+1)L\norm{\hat{u}}_\infty$, $s:=M_c(m+1)\bar{R}$, and $\bar{R}:=\max\{R,T\}$.
\endth

\begpr
Applying Theorem~\ref{th:GK-iterated-integral-approx-improved} gives the following:
\begin{align*}
\abs{F_c[u](T)-\hat{F}_c^J[\hat{u}](L)}
&\leq\sum_{j=0}^J\sum_{\eta\in X^j}\abs{(c,\eta)}\abs{E_\eta[u](T,0)-S_{\eta}[\hat{u}](L)}+ \\
&\hspace*{0.28in}\sum_{j=J+1}^\infty \sum_{\eta\in X^j} \abs{(c,\eta)}\abs{E_\eta[u](T,0)} \\
&\lesssim\sum_{j=2}^J K_cM_c^j(m+1)^j \, \frac{T^j}{L}\frac{\norm{\hat{u}/\Delta}_\infty^j}{2(j-2)!}+ \\
&\hspace*{0.28in}\sum_{j=J+1}^\infty K_cM_c^j(m+1)^j \,\frac{\bar{R}^j}{j!} \\
&=\frac{K_c}{2L}\sum_{j=0}^J (M_c(m+1)L\norm{\hat{u}}_\infty)^{j+2} \frac{1}{j!}+ \\
&\hspace*{0.28in}K_c\sum_{j=J+1}^\infty (M_c(m+1)\bar{R})^j\frac{1}{j!}\\
&=\frac{K_c}{2L}\expup^{\hat{s}}\hat{s}^2\Gamma(J+1,\hat{s})+K_c\expup^s(1-\Gamma(J+1,s)) \\
&=\hat{e}(J)+e(J),
\end{align*}
where the identity $\sum_{j=0}^J s^j/j!=\expup^s\Gamma(J+1,s)$ has been used \cite[Chapter 8.35]{Gradshteyn-Ryzhik_80}.
\endpr

Analogous to the local case, the error bound in the previous theorem can be simplified when $L\gg J\gg 1$.

\begco \label{co:DT-FO-appproximation-GC}
Let $c\in\allseriesellGC$ with growth constants $K_c,M_c>0$.
If $u\in B_1^m(R)[0,T]$ and
$L\gg J\gg 1$ then
\begdi
\abs{F_c[u](T)-\hat{F}_c^J[\hat{u}](L)}
\lesssim \frac{K_c}{2L}\expup^{M_c(m+1)L\norm{\hat{u}}_\infty}(M_c(m+1)L\norm{\hat{u}}_\infty)^2.
\enddi
\endco
\begpr
The upper bound follows directly from Theorem~\ref{th:DT-FO-appproximation-GC} using the fact that $\lim_{J\rightarrow+\infty}$ $\Gamma(J,s)=1$
(since $\Gamma(J+1,s)=\expup^{-s}\sum_{j=0}^J s^j/j!$, $J\geq 0$).
\endpr

\begin{figure}[t]
\begin{center}
\includegraphics[scale=0.5]{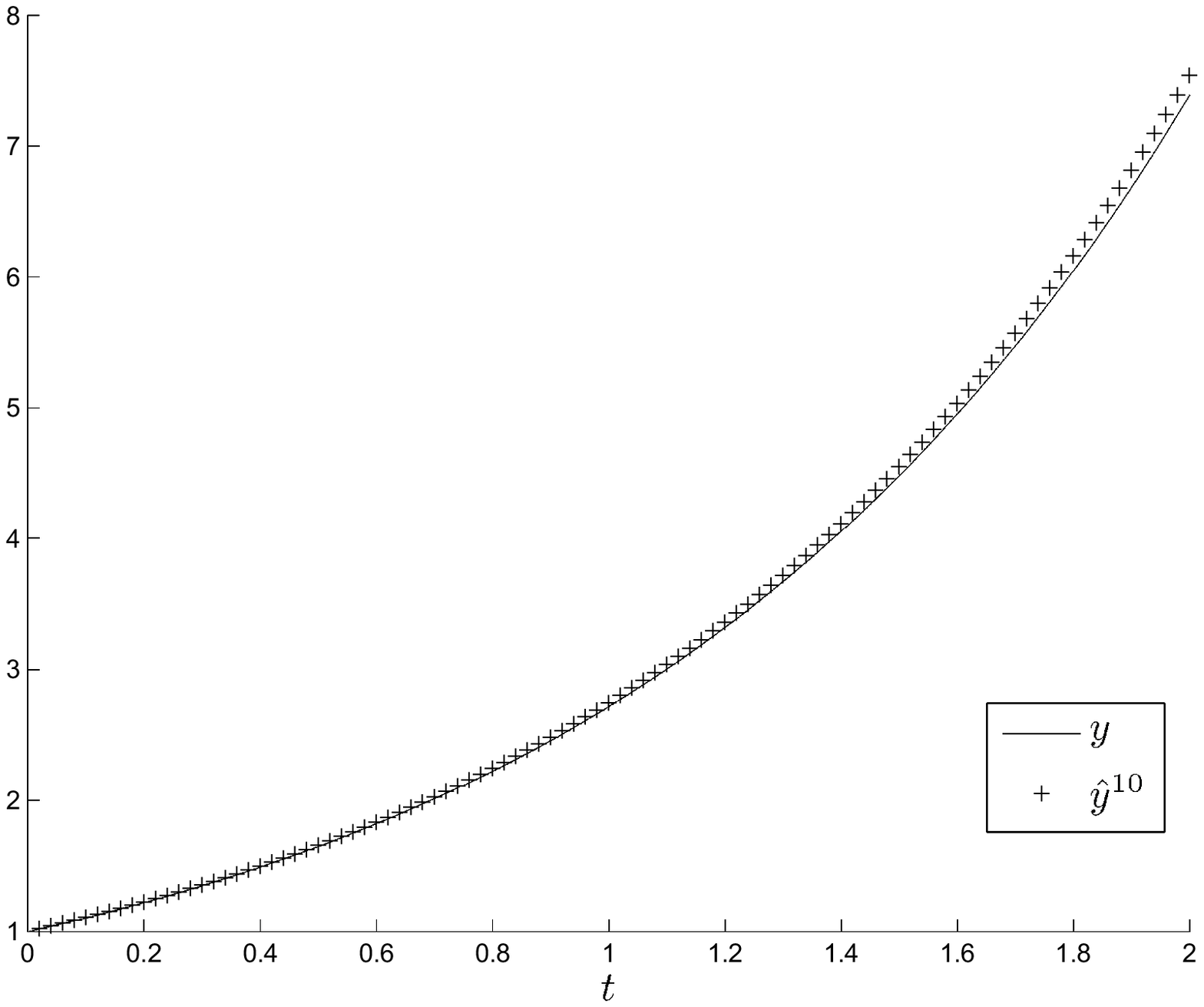}
\caption{Simulation comparing $y=F_c[1]$ to its approximation $\hat{y}^{10}=\hat{F}_c^{10}[\Delta]$ in
Example~\ref{ex:rational-Ferfera-example-DT-FO-approximation}, case 3.}
\label{fig:CTvsDT-rational-Ferfera-systems}
\end{center}
\end{figure}

\begin{figure}[t]
\begin{center}
\includegraphics[scale=0.5]{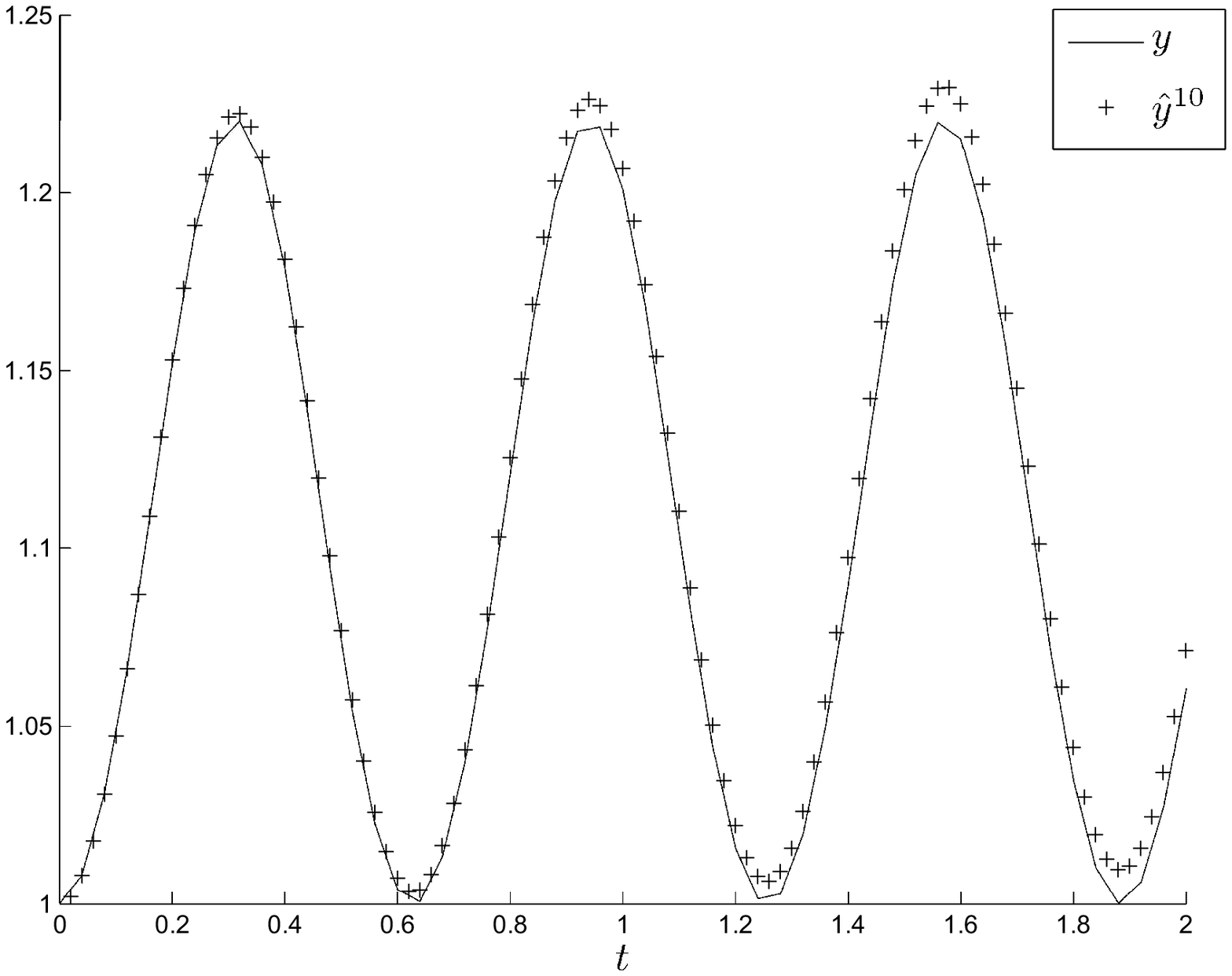}
\caption{Simulation comparing $y(t)=F_c[\sin(10t)]$ to its approximation $\hat{y}^{10}=\hat{F}_c^{10}[\hat{u}]$ in
Example~\ref{ex:rational-Ferfera-example-DT-FO-approximation}, case 6.}
\label{fig:CTvsDT-rational-Ferfera-systems2}
\end{center}
\end{figure}

\begex \label{ex:rational-Ferfera-example-DT-FO-approximation}
Consider the globally convergent series $c=\sum_{k\geq 0} x_1^k$ so that $K_c=M_c=1$.
In this case $F_c$ has the state space realization
\begeq \label{eq:CT-rational-Ferfera-system-realization}
\dot{z}=u,\;\; z(0)=0,\;\; y=\expup^z
\endeq
since
\begdi
y(t)=\sum_{j=0}^{\infty} (E_{x_1}[u](t))^j\frac{1}{j!}= \sum_{j=0}^{\infty} E_{x_1^{\shuffle j}\frac{1}{j!}}[u](t)
=\sum_{j=0}^{\infty}\,E_{x_1^j}[u](t) = F_c[u](t)
\enddi
for all $t\geq 0$. If the constant input $u=1$ is applied over the interval $[0,T]$ then $y(T)=\expup^{T}$.
The discrete-time approximation at $T=L\Delta$ is
\begin{align*}
\hat{y}^J(L)&=\hat{F}^J_c[\Delta](L)=\sum_{j=0}^J S_{x_1^j}[\Delta](L)
= \sum_{j=0}^J \Delta^j\sum_{k_1+k_2+\cdots +k_L=j} 1 =\sum_{j=0}^J \Delta^j {L+j-1 \choose j} \\
&=\sum_{j=0}^J \frac{\Delta^j}{j!} \left(L^j+\frac{(j-1)j}{2}L^{j-1}+\cdots+(j-1)!L\right) \\
&\approx \sum_{j=0}^J \frac{T^j}{j!} + \frac{1}{2L}\sum_{j=2}^J \frac{T^j}{(j-2)!} \\
&=[F_c[1](T)-e(J)]+\hat{e}(J),
\end{align*}
which is consistent with Theorem \ref{th:DT-FO-appproximation-GC} and again the worst case scenario in terms of approximating
the iterated integrals. The outputs $y$ and $\hat{y}$ were computed numerically over the interval $[0,2]$ for various choices of
$u$, $L$, and $J$.
This data is summarized in Table~\ref{tbl:simulation-summary-rational-case}, and the corresponding
plots for cases 3 and 6 are shown in
Figures~\ref{fig:CTvsDT-rational-Ferfera-systems} and
\ref{fig:CTvsDT-rational-Ferfera-systems2}, respectively.
As in the previous example, most of the error in
the approximation is due to the term $\hat{e}(J)$, and the constant input case yields an error that is
approximately upper bounded by $\hat{e}(J)+e(J)$. The error bound for the sinusoidal input is again conservative.
\endex

\section{Approximating Rational Operators}
\label{sec:rational-operators}

In the case where $F_c$ is a rational operator, it is shown in this section that the approximation $\hat{F}_c$
can be computed {\em without} the need for truncation. This is due exclusively to the fact that the generating series for
such an operator has structure which is not available in general, namely, a linear representation as
described in Definition~\ref{def:linear-representation}. The main idea is to use this representation to construct a
discrete-time state space realization for $\hat{F}_{c}$.
Later it will be shown that this technique is directly related to a specific discretization of the corresponding bilinear state space realization of
$F_c$. But the connection only becomes apparent in retrospect. For simplicity, the focus will be on the single-output case.
As motivation, consider the following simple example.
\begex
If $c=x_{i_3}x_{i_2}x_{i_1}\in X^\ast$ then the corresponding discrete-time Fliess operator is
$\hat{y}=S_{x_{i_3}x_{i_2}x_{i_1}}[\hat{u}]$.
Define the state $\hat{z}_1=S_{x_{i_1}}[\hat{u}]$ so that
\begin{displaymath}
\hat{z}_1(N+1)=\hat{z}_1(N)+\hat{u}_{i_1}(N+1).
\end{displaymath}
Similarly, if $\hat{z}_2=S_{x_{i_2i_1}}[\hat{u}]$ then
\begin{align*}
\hat{z}_2(N+1)&= \hat{z}_2(N)+\hat{u}_{i_2}(N+1)\hat{z}_1(N+1) \\
&= \hat{z}_2(N)+\hat{z}_1(N)\hat{u}_{i_2}(N+1)+
\hat{u}_{i_2}(N+1)\hat{u}_{i_1}(N+1).
\end{align*}
Finally, setting $\hat{y}=\hat{z}_3=S_{x_{i_3}x_{i_2}x_{i_1}}[\hat{u}]$ gives
\begin{align*}
\hat{z}_3(N+1)&=\hat{z}_3(N)+\hat{z}_2(N)\hat{u}_{i_3}(N+1)+
\hat{z}_1(N)\hat{u}_{i_3}(N+1)\hat{u}_{i_2}(N+1)+\\
&\hspace*{0.28in}\hat{u}_{i_3}(N+1)\hat{u}_{i_2}(N+1)\hat{u}_{i_1}(N+1).
\end{align*}
\endex

This triangular polynomial system is clearly not input-affine, as would be the case for the analogous
continuous-time input-output system $y=E_{x_{i_3}x_{i_2}x_{i_1}}[u]$, but the realization is {\em state affine}
in the following sense.

\begde \label{def:state-affine-polynomial-system}
A discrete-time state space realization is \bfem{polynomial input} and \bfem{state affine} if its transition map has the form
\begin{displaymath}
\hat{z}_i(N+1)=\sum_{j=1}^n p_{ij}(\hat{u}(N+1))\hat{z}_j(N)+q_i(\hat{u}(N+1)),
\end{displaymath}
$i=1,2,\ldots,n$, where $\hat{z}(N)\in\mathbb{R}^n$, $\hat{u}=[\hat{u}_0,\hat{u}_1,\ldots,\hat{u}_m]^T$, $p_{ij}$ and $q_i$ are polynomials,
and the output map $h:\hat{z}\mapsto \hat{y}$ is linear.
\endde

Polynomial input, state affine systems constitute an important class of discrete-time systems as first observed by Sontag in \cite[Chapter~V]{Sontag_79}.
The fact that $\hat{u}(N+1)$ appears in the transition map instead of $\hat{u}(N)$, as is more common, has no serious consequences here.
It will turn out, however, that if $c$ is rational instead of being merely polynomial, a more general class of state space realization
is required, one where rational functions of the input are admissible.

\begde \label{def:state-affine-rational-system}
A discrete-time state space realization is \bfem{rational input} and \bfem{state affine} if its transition map has the form
\begin{displaymath}
\hat{z}_i(N+1)=\sum_{j=1}^n r_{ij}(\hat{u}(N+1))\hat{z}_j(N)+s_i(\hat{u}(N+1)),
\end{displaymath}
$i=1,2,\ldots,n$, where $\hat{z}(N)\in\mathbb{R}^n$, $\hat{u}=[\hat{u}_0,\hat{u}_1,\ldots,\hat{u}_m]^T$,
$r_{ij}$ and $s_i$ are rational functions, and the output map $h:\hat{z}\mapsto \hat{y}$ is linear.
\endde

The main theorem of the section is below.

\begth \label{th:rational-state-affine-realization}
Let $c\in\mathbb{R}\langle\langle X \rangle\rangle$ be a rational series over $X=\{x_0,x_1,\ldots,x_m\}$ with
representation $(\mu,\gamma,\lambda)$.
Then $\hat{y}=\hat{F}_c[\hat{u}]$ has a finite dimensional rational input and state affine realization on $B_\infty^{m+1}[0,N_f](\hat{R})$ for
any $N_f> 0$ provided
$\hat{R}>0$ is sufficiently small.
\endth

Before giving the proof, some preliminary results are needed.

\begle \label{le:one-step-Fc-indentity}
For any $c\in\allseries$ it follows that
\begdi
\hat{F}_c[\hat{u}](N+1)=\hat{F}_c[\hat{u}](N)+\sum_{j=0}^m \hat{u}_j(N+1)\,\hat{F}_{x_j^{-1}(c)}[\hat{u}](N+1).
\enddi
\endle

\noindent{\em Proof:}
Observe that
\begin{align*}
\hat{F}_c[\hat{u}](N+1)
&=\sum_{\eta\in X^\ast} (c,\eta) S_\eta[\hat{u}](N+1) \\
&=\sum_{j=0}^m\sum_{\eta\in X^\ast} (c,x_j\eta) \sum_{k=0}^{N+1} \hat{u}_j(k)S_\eta[\hat{u}](k) \\
&=\sum_{j=0}^m\sum_{\eta\in X^\ast} (c,x_j\eta) \left[\sum_{k=0}^{N} \hat{u}_j(k)S_\eta[\hat{u}](k)+
\hat{u}_j(N+1)S_\eta[\hat{u}](N+1)\right] \\
&=\hat{F}_c[\hat{u}](N)+\sum_{j=0}^m \hat{u}_j(N+1)\sum_{\eta\in X^\ast} (x_j^{-1}(c),\eta)S_\eta[\hat{u}](N+1) \\
&=\hat{F}_c[\hat{u}](N)+\sum_{j=0}^m \hat{u}_j(N+1)\hat{F}_{x_j^{-1}(c)}[\hat{u}](N+1).
\end{align*}
\endpr

The next theorem hints at the well known dichotomy between time-reversible and non-time-reversible discrete-time systems.
That is, while every continuous-time state space realization can be run in reverse time, this is definitely not the case
for discrete-time systems. The system in the following theorem will only be time-reversible under certain conditions.

\begth
Let $c\in\mathbb{R}\langle\langle X \rangle\rangle$ be a rational series over $X=\{x_0,x_1,\ldots,x_m\}$.
Then $\hat{y}=\hat{F}_c[\hat{u}]$ has a finite dimensional backward-in-time bilinear realization for any
input sequence $\hat{u}$ defined over $[0,N_f]$.
\endth

\noindent{\em Proof:} Since $c$ is rational, it follows from Theorem~\ref{th:stable-rational-subspace}
that a stable $n$ dimension subspace $V$ of $\allseries$ exists which contains $c$. Let $\bar{c}_k$, $k=1,2,\ldots,n$ be a basis for $V$ so that
$c=\sum_{k=1}^n\lambda_k \bar{c}_k$ with $\lambda_k\in\re$.
Furthermore,
for any $x_j\in X$ it follows that
\begdi
x_j^{-1}(\bar{c}_k)=\sum_{l=1}^n \mu_{kl}(x_j)\, \bar{c}_l,
\enddi
where $\mu_{kl}(x_j)\in\re$. Define the state variables
$\bar{z}_k(N)=\hat{F}_{\bar{c}_k}[\hat{u}](N_f-N)$, $k=1,2,\ldots,n$
for $N\in [0,N_f]$. Then
\begdi
\hat{y}(N)=\hat{F}_c[\hat{u}](N)=\sum_{k=1}^n \lambda_k \hat{F}_{\bar{c}_k}[\hat{u}](N)=\sum_{k=1}^n \lambda_k \bar{z}_k(N_f-N)
\enddi
and
\begdi
\bar{z}_k(N_f)=\hat{F}_{\bar{c}_k}[\hat{u}](0)=(\bar{c}_k,\emptyset)S_{\emptyset}[\hat{u}](0)=(\bar{c}_k,\emptyset)=:\gamma_k.
\enddi
Now from Lemma~\ref{le:one-step-Fc-indentity}
\begin{align*}
\bar{z}_k(N_f-N-1)
&=\bar{z}_k(N_f-N)\!+\!\sum_{j=0}^m \hat{u}_j(N+1)\,\hat{F}_{x_j^{-1}(\bar{c}_k)}[\hat{u}](N+1) \\
&=\bar{z}_k(N_f-N)\!+\!\sum_{j=0}^m \hat{u}_j(N+1)\sum_{l=1}^n \mu_{kl}(x_j)\hat{F}_{\bar{c}_l}[\hat{u}](N+1) \\
&=\bar{z}_k(N_f-N)\!+\!\sum_{j=0}^m \hat{u}_j(N+1)\sum_{l=1}^n \mu_{kl}(x_j)\bar{z}_l(N_f-N-1) \\
&=\bar{z}_k(N_f-N)\!+\!\sum_{j=0}^m \hat{u}_j(N+1)\,[A_j\bar{z}(N_f-N-1)]_k,
\end{align*}
where $A_j\in\re^{n\times n}$, $j=0,1,\ldots,m$ has components $[A_j]_{kl}=\mu_{kl}(x_j)$,
and $\bar{z}$ is the column vector with $\bar{z}_k$ as its $k$-th component.
Therefore, for $N=N_f-1, N_f-2,\ldots,0$ it follows that
\begdi
\bar{z}(N_f-N)=\left[I-\sum_{j=0}^m A_j\hat{u}_j(N+1)\right]\bar{z}(N_f-N-1),
\enddi
with $\bar{z}(N_f)=\gamma$ and $\hat{y}(N)=\lambda \bar{z}(N_f-N)$, or equivalently, setting $\hat{z}(N)=\bar{z}(N_f-N)$ gives for $N=N_f-1,N_f-2,\ldots,0$
\begeq
\hat{z}(N)=\left[I-\sum_{j=0}^m A_j \hat{u}_j(N+1)\right]\hat{z}(N+1) \label{eq:backwards-bilinear-system}
\endeq
with $\hat{z}(0)=\gamma$ and $\hat{y}(N)=\lambda \hat{z}(N)$ as claimed.
{\hfill\bull}

\vspace*{0.1in}

The proof of the main result follows from introducing conditions on $\hat{u}$ so that
system~\rref{eq:backwards-bilinear-system}  is time-reversible.
Bilinearity is lost in the process, but the forward-in-time system is rational input and state affine.

\vspace*{0.1in}

\noindent{\em Proof of Theorem~\ref{th:rational-state-affine-realization}:}
If $\hat{u}\in B_\infty^{m+1}[0,N_f](\hat{R})$, and $\hat{R}$ is
sufficiently small, then the transition matrix $I-\sum_{j=0}^mA_j \hat{u}_j(N+1)$ of system
\rref{eq:backwards-bilinear-system} is nonsingular.
In which case, the forward-in-time system
\begdi
\hat{z}(N+1)=\left[I-\sum_{j=0}^mA_j \hat{u}_j(N+1)\right]^{-1}\hat{z}(N)
\enddi
is well defined over $[0,N_f]$ and clearly state affine and rational in $\hat{u}$. Furthermore, by design
$\hat{y}=F_c[\hat{u}]=\lambda \hat{z}$ over the interval $[0,N_f]$.
{\hfill\bull}

\begex \label{ex:rational-Ferfera-example-DT-FO-approximation-revisited}
Reconsider the rational Fliess operator in Example~\ref{ex:rational-Ferfera-example-DT-FO-approximation} where $c=\sum_{k\geq 0} x_1^k$.
Clearly, $x_0^{-1}(c)=0$, $x_1^{-1}(c)=c$, and $(c,\emptyset)=1=1\cdot 1=\lambda\gamma$.
Thus, $\hat{F}_c$ has the $n=1$ dimensional
rational and state affine
realization
\begin{equation}
\label{eq:DT-rational-Ferfera-system-realization}
\hat{z}(N+1)=(1-\hat{u}(N+1))^{-1}\hat{z}(N),\;\; \hat{z}(0)=1,\;\;
\hat{y}(N)=\hat{z}(N)
\end{equation}
provided $\norm{\hat{u}}_{\infty}< 1$.
Since
\begdi
\hat{z}(N+1)=\sum_{i=0}^{\infty} \hat{u}^i(N+1)\hat{z}(N),
\enddi
if follows for $N\geq 0$ that
\begdi
\hat{y}(N)=\prod_{k=1}^N (1-\hat{u}(k))^{-1}\hat{z}(0)
=\sum_{i_1,\ldots,i_N=0}^\infty \hat{u}^{i_N}(N)\hat{u}^{i_{N-1}}(N-1)\cdots \hat{u}^{i_1}(1),
\enddi
where the product is defined to be unity when $N=0$.
For example,
\begin{align*}
\hat{y}(0)&=1 \\
\hat{y}(1)&=1+\hat{u}(1)+\hat{u}^2(1)+\hat{u}^{3}(1)+\cdots \\
\hat{y}(2)&=(1+\hat{u}(2)+\hat{u}^2(2)+\cdots)(1+\hat{u}(1)+\hat{u}^2(1)+\cdots) \\
&=1+(\hat{u}(1)+\hat{u}(2))+(\hat{u}^2(1)+\hat{u}(2)\hat{u}(1)+\hat{u}^2(2))+ \\
&\hspace*{0.2in}
(\hat{u}^3(1)+\hat{u}^2(2)\hat{u}(1)+\hat{u}(2)\hat{u}^2(1)+\hat{u}^3(2))+\cdots \\
&\hspace*{0.1in}\vdots
\end{align*}
This solution can be checked independently by simply applying the definition of $\hat{F}_{c}$. That is,
\begdi
\hat{y}(N)=S_{\emptyset}[\hat{u}](N)+S_{x_1}[\hat{u}](N)+S_{x_1^2}[\hat{u}](N)+\cdots,
\enddi
so that
\begin{align*}
\hat{y}(0)&=1 \\
\hat{y}(1)&=1+\hat{u}(1)+\hat{u}^2(1)+\hat{u}^{3}(1)+\cdots \\
\hat{y}(2)&=1+(\hat{u}(1)+\hat{u}(2))+(\hat{u}^2(1)+\hat{u}(2)\hat{u}(1)+\hat{u}^2(2))+ \\
&\hspace*{0.2in}
(\hat{u}^3(1)+\hat{u}^2(2)\hat{u}(1)+\hat{u}(2)\hat{u}^2(1)+\hat{u}^3(2))+\cdots \\
&\hspace*{0.1in}\vdots
\end{align*}
Not surprisingly, the plots of $\hat{y}$ generated from system~\rref{eq:DT-rational-Ferfera-system-realization} are indistinguishable from those shown in
Figures~\ref{fig:CTvsDT-rational-Ferfera-systems} and \ref{fig:CTvsDT-rational-Ferfera-systems2},
which were generated directly from the definition of $\hat{F}_c^J$. It also should be noted
that $F_c$, being rational, has a bilinear realization
\begdi
\dot{\tilde{z}}=\tilde{z}u,\;\; \tilde{z}(0)=1,\;\; y=\tilde{z},
\enddi
which is related to the realization \rref{eq:CT-rational-Ferfera-system-realization} by the coordinate transformation
$\tilde{z}=\expup^z$.
For small $\Delta>0$ observe
\begin{align*}
\tilde{z}((N+1)\Delta)&=\tilde{z}(N\Delta)+\int_{N\Delta}^{(N+1)\Delta} \tilde{z}(t)u(t)\,dt \\
&\approx \tilde{z}(N\Delta)+\int_{N\Delta}^{(N+1)\Delta} \hspace*{-0.1in}u(t)\,dt\,\, \tilde{z}((N+1)\Delta)\\
&=\tilde{z}(N\Delta)+\hat{u}(N+1)\,\tilde{z}((N+1)\Delta),
\end{align*}
and therefore, letting $\hat{z}(N)=\tilde{z}(N\Delta)$, this particular discretized system
\begdi
\hat{z}(N+1)=(1-\hat{u}(N+1))^{-1}\hat{z}(N)
\enddi
has the form of \rref{eq:DT-rational-Ferfera-system-realization}.
\endex

\begex
The previous example can be generalized by noting that
\begin{align*}
\hat{z}(N+1)
&=\left[I-\sum_{j=0}^mA_j \hat{u}_j(N+1)\right]^{-1}\hat{z}(N) \\
&=\left[\sum_{k=0}^\infty\sum_{j_0,\ldots,j_k=0}^m A_{j_0}A_{j_1}\cdots A_{j_k} \hat{u}_{j_0}(N+1)
\hat{u}_{j_1}(N+1)\cdots \hat{u}_{j_k}(N+1)\right] \hat{z}(N)\\
&=:\sum_{\eta=x_{j_0}\cdots x_{j_k}\in X^\ast} A_{\eta}\hat{u}_\eta(N+1) \hat{z}(N).
\end{align*}
In which case,
\begin{align*}
\hat{y}(N)&=\sum_{\eta_N,\ldots,\eta_1\in X^\ast} \lambda A_{\eta_N}\cdots A_{\eta_1}\gamma\; \hat{u}_{\eta_N}(N+1)\cdots \hat{u}_{\eta_1}(1) \\
&=\sum_{\eta_N,\ldots,\eta_1\in X^\ast} (c,\eta_N\cdots \eta_1) \hat{u}_{\eta_N}(N+1)\cdots \hat{u}_{\eta_1}(1).
\end{align*}
This form of the discrete-time input-output map comes from a specific discretization of the underlying continuous-time realization
\begdi
\dot{z}(t)=\sum_{j=0}^m A_j z(t)u_j, \;\;z(0)=\gamma,\;\;
y(t)=\lambda z(t),
\enddi
namely,
\begdi
z((N+1)\Delta)\approx z(N\Delta )+\sum_{j=0}^m A_j\hat{u}_j(N+1)\,z((N+1)\Delta)
\enddi
so that
\begdi
\hat{z}(N+1)=\left[I-\sum_{j=0}^m A_j\hat{u}_j(N+1)\right]^{-1}\hat{z}(N).
\enddi
\endex

\section{Conclusions}

This paper described how to approximate Fliess operators with iterated sums and gave explicit
achievable error bounds for the locally and globally convergent cases.
For the special case of
rational Fliess operators, it was shown that the method can be realized via a rational input and state affine discrete-time state space model.
This model avoids the truncation error and
can also be derived from a specific discretization of a continuous-time bilinear realization of the rational Fliess operator.

\newpage

\begin{table}[h]
\vspace*{4.0in}
\hspace*{0.3in}
\rotatebox{90}{
\begin{minipage}{2.8in}
\begin{center}
\caption{Summary of simulation results for Example~\ref{ex:non-rational-Ferfera-DT-FO-approximation}.}
\label{tbl:simulation-summary-nonrational-case}
\renewcommand{\arraystretch}{1.5}
{\small
\begin{tabular}{|c|c|c|c|c|c|c|c|c|c|c|c|c|c|} \hline
\multicolumn{1}{|c|}{case} &
\multicolumn{1}{|c|}{$u$} &
\multicolumn{1}{|c|}{$T$} &
\multicolumn{1}{|c|}{$L$} &
\multicolumn{1}{|c|}{$\Delta$} &
\multicolumn{1}{|c|}{$J$} &
\multicolumn{1}{|c|}{$\norm{\hat{u}}_\infty$} &
\multicolumn{1}{|c|}{$s$} &
\multicolumn{1}{|c|}{$\hat{s}$} &
\multicolumn{1}{|c|}{$y(T)$} &
\multicolumn{1}{|c|}{$\hat{y}^J(L)$} &
\multicolumn{1}{|c|}{$\hat{y}^J(L)-y(T)$} &
\multicolumn{1}{|c|}{$\hat{e}(J)$} &
\multicolumn{1}{|c|}{$e(J)$} \\ \hline\hline
1 & $1$ & 0.5 &50 & 0.0100 & 10 & 0.0100 & 0.5000 & 0.5000 & 2.0000 & 2.0412 & 0.0412          & 0.0355 & 9.7656$\times10^{-4}$  \\ \hline
2 & $1$ & 0.5 &50 & 0.0100 & 20 & 0.0100 & 0.5000 & 0.5000 & 2.0000 & 2.0448 & 0.0448          & 0.0400 & 9.5367$\times10^{-7}$  \\ \hline
3 & $1$ & 0.5 &100 & 0.0050 & 10 & 0.0050 & 0.5000 & 0.5000 & 2.0000 & 2.0192 & 0.0192         & 0.0177 & 9.7656$\times10^{-4}$  \\ \hline
4 & $\sin(20t)$ & 0.5 &50 & 0.0100 & 10 & 0.0099 & 0.5000 & 0.4975 & 1.1009 & 1.1041 & 0.0032  & 0.0347 & 9.7656$\times10^{-4}$  \\ \hline
5 & $\sin(20t)$ & 0.5 &50 & 0.0100 & 20 & 0.0099 & 0.5000 & 0.4975 & 1.1009 & 1.1041 & 0.0032  & 0.0390 & 9.5367$\times10^{-7}$  \\ \hline
6 & $\sin(20t)$ & 0.5 &100 & 0.0050 & 10 & 0.0050 & 0.5000 & 0.4994 & 1.1011 & 1.1028 & 0.0017 & 0.0176 & 9.7656$\times10^{-4}$  \\ \hline
\end{tabular}
}
\end{center}
\end{minipage}
}
\end{table}

\begin{table}[h]
\vspace*{-3.03in}
\hspace*{2.4in}
\rotatebox{90}{
\begin{minipage}{2.8in}
\begin{center}
\caption{Summary of simulation results for Example~\ref{ex:rational-Ferfera-example-DT-FO-approximation}.}
\label{tbl:simulation-summary-rational-case}
\renewcommand{\arraystretch}{1.5}
{\small
\begin{tabular}{|c|c|c|c|c|c|c|c|c|c|c|c|c|c|} \hline
\multicolumn{1}{|c|}{case} &
\multicolumn{1}{|c|}{$u$} &
\multicolumn{1}{|c|}{$T$} &
\multicolumn{1}{|c|}{$L$} &
\multicolumn{1}{|c|}{$\Delta$} &
\multicolumn{1}{|c|}{$J$} &
\multicolumn{1}{|c|}{$\norm{\hat{u}}_\infty$} &
\multicolumn{1}{|c|}{$s$} &
\multicolumn{1}{|c|}{$\hat{s}$} &
\multicolumn{1}{|c|}{$y(T)$} &
\multicolumn{1}{|c|}{$\hat{y}^J(L)$} &
\multicolumn{1}{|c|}{$\hat{y}^J(L)-y(T)$} &
\multicolumn{1}{|c|}{$\hat{e}(J)$} &
\multicolumn{1}{|c|}{$e(J)$} \\ \hline\hline
1 & $1$ & 2 &50 & 0.0400 & 10 & 0.0400 & 2.0000 & 2.0000 & 7.3891 & 7.6989 & 0.3098          & 0.2956 & 6.1390$\times10^{-5}$   \\ \hline
2 & $1$ & 2 &50 & 0.0400 & 20 & 0.0400 & 2.0000 & 2.0000 & 7.3891 & 7.6991 & 0.3100          & 0.2956 & 4.5119$\times10^{-14}$  \\ \hline
3 & $1$ & 2 &100 & 0.0200 & 10 & 0.0200 & 2.0000 & 2.0000 & 7.3891 & 7.5403 & 0.1512         & 0.1478 & 6.1390$\times10^{-5}$   \\ \hline
4 & $\sin(10t)$ & 2 &50 & 0.0400 & 10 & 0.0392 & 2.0000 & 1.9601 & 1.0601 & 1.0803 & 0.0202  & 0.2728 & 6.1390$\times10^{-5}$   \\ \hline
5 & $\sin(10t)$ & 2 &50 & 0.0400 & 20 & 0.0392 & 2.0000 & 1.9601 & 1.0601 & 1.0803 & 0.0202  & 0.2728 & 4.5119$\times10^{-14}$  \\ \hline
6 & $\sin(10t)$ & 2 &100 & 0.0200 & 10 & 0.0199 & 2.0000 & 1.9899 & 1.0607 & 1.0711 & 0.0104 & 0.1448 & 6.1390$\times10^{-5}$   \\ \hline
\end{tabular}
}
\end{center}
\end{minipage}
}
\end{table}

\end{document}